# Hankel determinants of convolution powers of Catalan numbers revisited

Johann Cigler

**Abstract**

Using a slightly generalized result of George Andrews and Jet Wimp this note gives a simple computational proof of some Hankel determinants of backwards shifts of convolution powers of Catalan numbers and obtains analogous results for Narayana polynomials.

## 1. Introduction

In [5] and [6] I obtained conjectures for Hankel determinants of backwards shifts of convolution powers of Catalan numbers. Recently Markus Fulmek [8] gave bijective proofs of these conjectures. In the present note I use the following Lemma ([3], Proposition 2.5)), which slightly generalizes a result of George Andrews and Jet Wimp [2], for another proof of my Conjectures and for analogous results for Narayana polynomials.

**Lemma**

Let $s(x) = \sum_{n \geq 0} s_n x^n$ with $s_0 = 1$ and $t(x) = \dfrac{1}{s(x)} = \sum_{n \geq 0} t_n x^n$.

Setting $s_n = t_n = 0$ for $n < 0$ we get for $M \in \mathbb{N}$

$$\det\left(s_{i+j-M}\right)_{i,j=0}^{N+M} = (-1)^{N+\binom{M+1}{2}} \det\left(t_{i+j+M+2}\right)_{i,j=0}^{N-1}. \tag{1}$$

## 2. Hankel determinants for convolutions of Catalan numbers

Let

$$c(x) = \frac{1-\sqrt{1-4x}}{2x} = \sum_{n \geq 0} C_n x^n \tag{2}$$

be the generating function of the Catalan numbers

$$C_n = \frac{1}{n+1}\binom{2n}{n}. \tag{3}$$

It is well known that

$$c(x)^k = \sum_{n \geq 0} C_{k,n} x^n = \sum_{n \geq 0} \frac{k}{n+k}\binom{2n+k-1}{n} x^n \tag{4}$$

and that the convolution powers $C_{k,n}$ of the Catalan numbers can be interpreted as the number of all non-negative lattice paths in $\mathbb{Z}^2$ with up-steps $U = (1,1)$ and down-steps $D = (1,-1)$, which we will call non-negative up-down paths, from $(0,0)$ to $(2n+k-1, k-1)$.

To verify this fact, observe that

$$1 + xc(x)^2 = 1 + \frac{\left(1 - \sqrt{1-4x}\right)^2}{4x} = 1 + \frac{2 - 4x - 2\sqrt{1-4x}}{4x} = \frac{1 - \sqrt{1-4x}}{2x} = c(x) \quad (5)$$

implies

$$x^k c\left(x^2\right)^{k+1} = x\left(x^{k-1} c\left(x^2\right)^k\right) + x\left(x^{k+1} c\left(x^2\right)^{k+2}\right). \quad (6)$$

If we write $x^k c\left(x^2\right)^{k+1} = \sum_{j \geq 0} a(j,k) x^j$ then (6) is equivalent with

$$a(j,k) = a(j-1, k-1) + a(j-1, k+1). \quad (7)$$

The identity (7) shows that $a(j,k)$ can be interpreted as the number of all non-negative up-down paths from $(0,0)$ to $(j,k)$.

The series $c(x)^k$ satisfies

$$\left(xc(x)\right)^k + \frac{1}{c(x)^k} = L_k(1, -x), \quad (8)$$

where $L_n(x,s)$ are the bivariate Lucas polynomials which are defined by

$$\begin{aligned} L_n(x,s) &= x L_{n-1}(x,s) + s L_{n-2}(x,s), \\ L_0(x,s) &= 2, \quad L_1(x,s) = x. \end{aligned} \quad (9)$$

As is well known, the Lucas polynomials satisfy

$$L_n\left(x+y, -xy\right) = x^n + y^n. \quad (10)$$

To show this it suffices to verify (10) for $n=0$ and $n=1$ and observe that
$x^n + y^n = (x+y)\left(x^{n-1} + y^{n-1}\right) - xy\left(x^{n-2} + y^{n-2}\right).$

To prove (8) let

$$h_k(x) = \left(xc(x)\right)^k + \frac{1}{c(x)^k}. \quad (11)$$

From (5) we get

$$h_1(x) = xc(x) + \frac{1}{c(x)} = 1. \quad (12)$$

(10) implies $h_k(x) = L_k\left(xc(x) + \frac{1}{c(x)}, -xc(x) \cdot \frac{1}{c(x)}\right) = L_k(1, -x).$

Let us note that $\deg L_k(1,-x) = \left\lfloor \dfrac{k}{2} \right\rfloor$. The first polynomials $L_k(1,-x)$ are $L_1(1,-x) = 1$, $L_2(1,-x) = 1-2x$, $L_3(1,-x) = 1-3x$, $L_4(1,-x) = 1-4x+2x^2$, $L_5(1,-x) = 1-5x+5x^2$, $L_6(1,-x) = 1-6x+9x^2-2x^3$.

For later use observe that

$$h_{2k}(x) = (xc(x))^{2k} + \frac{1}{c(x)^{2k}} = L_k\left((xc(x))^2 + \frac{1}{c(x)^2}, -x^2\right) = L_k\left(h_2(x), -x^2\right)$$

which implies

$$L_{2k}(1,-x) = L_k\left(1-2x, -x^2\right). \tag{13}$$

We extend $C_{k,n}$ to negative integers by setting $C_{k,n} = 0$ for $n < 0$ and define

$$D_{K,M}(N) := \det\left(C_{K,i+j+M}\right)_{i,j=0}^{N-1} \tag{14}$$

for $M \in \mathbb{Z}$ and integers $K, N \geq 1$. We set $D_{K,M}(0) = 1$ by definition.

The following Theorems have been conjectured in [5] and [6] and bijective proofs have been given by Markus Fulmek [8].

**Theorem 1.**

*For $m \in \mathbb{N}$ and a positive integer $k$ we have*

$$D_{2k,1-k-m}(N) = 0 \text{ for } N = 1, 2, \cdots, m+k-1, \tag{15}$$

*and for all $n \in \mathbb{N}$*

$$D_{2k,1-k-m}(n+m+k) = (-1)^{\binom{m+k}{2}} D_{2k,1-k+m}(n). \tag{16}$$

For example for $k = 2$ and $m = 1$ we get

$$\left(D_{4,-2}(n)\right)_{n\geq 0} = (1, 0, 0, -1, -1, 2, 2, -3, -3, 4, 4, -5, -5, \cdots)$$
$$\left(D_{4,0}(n)\right)_{n\geq 0} = (1, 1, -2, -2, 3, 3, -4, -4, 5, 5, -6, -6, \cdots).$$

**Theorem 2.**

*For $m \in \mathbb{N}$ and a positive integer $k$ we have*

$$D_{2k-1,2-k-m}(N) = 0 \text{ for } N = 1, 2, \cdots, m+k-2, \tag{17}$$

*and for all $n \in \mathbb{N}$*

$$D_{2k-1,2-k-m}(n+m+k-1) = (-1)^{\binom{m+k-1}{2}} D_{2k-1,1-k+m}(n). \tag{18}$$

For example, for $k = 2$, $m = 1$ we get

$$(D_{3,-1}(n))_{n\geq 0} = (1,0,-1,-1,0,1,1,0,-1,-1,0,1,1,\cdots),$$
$$(D_{3,0}(n))_{n\geq 0} = (1,1,0,-1,-1,0,1,1,0,-1,-1,0,1,\cdots).$$

**Proof.**

The identities (15) and (17) are trivially true because the first rows of the Hankel matrices vanish. For (16) and (18) we can apply the Lemma.

To prove (16) choose $s(x) = c(x)^{2k}$ and $M = m+k-1$, $N = n$. Then we get

$$\det\left(s_{i+j-M}\right)_{i,j=0}^{N+M} = \det\left(s_{i+j+1-m-k}\right)_{i,j=0}^{n+m+k-1} = D_{2k,1-k-m}(n+m+k).$$

By (8) we get $t(x) = L_{2k}(1,-x) - (xc(x))^{2k}$ and therefore $t_{2k+r} = -C_{2k,r}$ for $r \geq 0$. Note that for $k < n < 2k$ $t_n = 0$ and $C_{2k,n-2k} = 0$. Therefore $t_{r+m+k+1} = -C_{2k,r+m-k+1}$ for $r+m+1 > 0$.

This gives $\det\left(t_{i+j+M+2}\right)_{i,j=0}^{N-1} = \det\left(t_{i+j+m+k+1}\right)_{i,j=0}^{n-1} = \det\left(-C_{2k,m-k+1}\right)_{i,j=0}^{n-1} = (-1)^n D_{2k,m-k+1}(n).$

Thus (1) gives (16).

For Theorem 2 we choose $s(x) = c(x)^{2k-1}$ and $M = m+k-2$. Then we get

$$\det\left(s_{i+j-M}\right)_{i,j=0}^{N+M} = \det\left(s_{i+j+2-m-k}\right)_{i,j=0}^{n+m+k-2} = D_{2k-1,2-k-m}(n+m+k-1).$$

On the other hand, we get $t(x) = L_{2k-1}(1,-x) - (xc(x))^{2k-1}$ and therefore $t_{2k-1+r} = -C_{2k-1,r}$ and $t_r = 0$ for $k-1 < r < 2k-1$.

This gives $\det\left(t_{i+j+M+2}\right)_{i,j=0}^{N-1} = \det\left(t_{i+j+m+k}\right)_{i,j=0}^{n-1} = \det\left(-C_{2k-1,m-k+1}\right)_{i,j=0}^{n-1} = (-1)^n D_{2k-1,m-k+1}(n).$

**Remark**

As special cases we get the following analogs of the well-known identities

$$\det\left(C_{i+j}\right)_{i,j=0}^{n-1} = 1,$$
$$\det\left(C_{2,i+j}\right)_{i,j=0}^{n-1} = \det\left(C_{i+j+1}\right)_{i,j=0}^{n-1} = 1. \tag{19}$$

For $m=0$ (16) gives $D_{2k,1-k}(n+k) = (-1)^{\binom{k}{2}} D_{2k,1-k}(n)$ and therefore

$$D_{2k,1-k}(kn) = (-1)^{n\binom{k}{2}} \tag{20}$$

and $D_{2k,1-k}(n) = 0$ else, which has been obtained with other means in [7], Cor. 15.

(18) gives for $m=0$ $D_{2k-1,2-k}(n+k-1) = (-1)^{\binom{k-1}{2}} D_{2k-1,1-k}(n)$ and for $m=1$

$D_{2k-1,1-k}(n+k) = (-1)^{\binom{k}{2}} D_{2k-1,2-k}(n)$. By changing $k \to k+1$ we get for $k \geq 1$

$$D_{2k+1,-k}((2k+1)n) = (-1)^{kn},$$
$$D_{2k+1,-k}((2k+1)n+k+1) = (-1)^{kn+\binom{k+1}{2}}, \tag{21}$$
$$D_{2k+1,-k}(n) = 0 \quad \text{else}$$

and

$$D_{2k+1,1-k}((2k+1)n) = (-1)^{kn},$$
$$D_{2k+1,1-k}((2k+1)n+k) = (-1)^{kn+\binom{k}{2}}, \tag{22}$$
$$D_{2k+1,1-k}(n) = 0 \quad \text{else,}$$

which previously has been proved with other means in [7], Cor.12 and Cor. 13.

## 3. Hankel determinants for convolutions of Narayana polynomials

Let $C_n(t) = \sum_{k=0}^{n-1} \binom{n}{k}\binom{n-1}{k} \frac{1}{k+1} t^k$ for $n > 0$ and $C_0(t) = 1$ be the Narayana polynomials

$$(C_n(t))_{n \geq 0} = \left(1, 1, 1+t, 1+3t+t^2, 1+6t+6t^2+t^3, 1+10t+20t^2+10t^3+t^4, \cdots\right)$$

which for $t = 1$ reduce to the Catalan numbers.

By [9], (2.6) we get the generating function

$$c_0(x,t) = \sum_{n \geq 0} C_n(t) x^n = \frac{1+(t-1)x - \sqrt{1-2x(t+1)+x^2(t-1)^2}}{2xt}. \tag{23}$$

We will also consider another generating function

$$c_1(x,t) = 1 + t\sum_{n \geq 1} C_n(t) x^n = \frac{1-(t-1)x - \sqrt{1-2x(t+1)+x^2(t-1)^2}}{2x} \tag{24}$$

which satisfies

$$c_1(x,t) = 1 - t + t c_0(x,t). \tag{25}$$

Note that $c_0(x,1) = c_1(x,1) = c(x)$. As analogs of (5) we get

$$c_1(x,t) = 1 + xt c_0(x,t) c_1(x,t),$$
$$c_0(x,t) = 1 + x c_0(x,t) c_1(x,t), \tag{26}$$

because

$$4x^2tc_0(x,t)c_1(x,t) = \left(1-\sqrt{1-2x(t+1)+x^2(t-1)^2}\right)^2 - (t-1)^2x^2$$

$$= 1 - 2\sqrt{1-2x(t+1)+x^2(t-1)^2} + 1 - 2x(t+1) = 2\left(1-x(t+1)-\sqrt{1-2x(t+1)+x^2(t-1)^2}\right)$$

$$= -4x + 2\left(1-x(t-1)-\sqrt{1-2x(t+1)+x^2(t-1)^2}\right) = -4x + 4xc_1(x,t)$$

Let us consider the series

$$\begin{aligned} c(x,t)^{(2k)} &= c(x,t)^{(2k-1)}c_1(x,t), \\ c(x,t)^{(2k+1)} &= c(x,t)^{(2k)}c_0(x,t), \\ c(x,t)^{(0)} &= 1 \end{aligned} \tag{27}$$

and the polynomials

$$C_{k,n}(t) = [x^n]c(x,t)^{(k)}, \tag{28}$$

which will be called convolution powers of the Narayana polynomials. They have also been studied in [4] with a somewhat different notation.

By (27) we get

$$c(x,t)^{(k)} = c(x,t)^{(k-2)} c(x,t)^{(2)}. \tag{29}$$

This implies

$$xc(x,t)^{(2k)} = c(x,t)^{(2k-2)} xc(x,t)^{(2)} = c(x,t)^{(2k-2)}\left(c_0(x,t)-1\right)$$
$$= c(x,t)^{(2k-1)} - c(x,t)^{(2k-2)}$$

or

$$c(x,t)^{(2k-1)} = c(x,t)^{(2k-2)} + xc(x,t)^{(2k)}. \tag{30}$$

and

$$xtc(x,t)^{(2k+1)} = c(x,t)^{(2k-1)} xtc(x,t)^{(2)} = c(x,t)^{(2k-1)}\left(c_1(x,t)-1\right)$$
$$= c(x,t)^{(2k)} - c(x,t)^{(2k-1)}$$

or

$$c(x,t)^{(2k)} = c(x,t)^{(2k-1)} + xtc(x,t)^{(2k+1)}. \tag{31}$$

Thus

$$\begin{aligned} C_{2k,n}(t) &= C_{2k-1,n}(t) + tC_{2k+1,n-1}(t), \\ C_{2k+1,n}(t) &= C_{2k,n}(t) + C_{2k+2,n-1}(t). \end{aligned} \tag{32}$$

Of course, we have $C_{1,n}(t) = C_n(t)$. And by (26) we get

$$C_{2,n}(t) = C_{n+1}(t). \tag{33}$$

For example,

$$\left(C_{2,n}(t)\right)_{n\geq 0} = \left(1, 1+t, 1+3t+t^2, 1+6t+6t^2+t^3, 1+10t+20t^2+10t^3+t^4, \cdots\right),$$

$$\left(C_{3,n}(t)\right)_{n\geq 0} = \left(1, 2+t, 3+5t+t^2, 4+14t+9t^2+t^3, 5+30t+40t^2+14t^3+t^4, \cdots\right).$$

Let $x^k c(x^2,t)^{(k+1)} = \sum_{j\geq 0} a(j,k,t) x^j$. Then (30) and (31) imply

$$x^{2k} c(x^2,t)^{(2k+1)} = x\left(x^{2k-1} c(x^2,t)^{(2k)}\right) + x\left(x^{2k+1} c(x^2,t)^{(2k+2)}\right)$$

which gives

$$a(n, 2k, t) = a(n-1, 2k-1, t) + a(n-1, 2k+1, t) \tag{34}$$

and

$$x^{2k-1} c(x^2,t)^{(2k)} = x\left(x^{2k-2} c(x^2,t)^{(2k-1)}\right) + xt\left(x^{2k} c(x^2,t)^{(2k+1)}\right)$$

which gives

$$a(n, 2k-1, t) = a(n-1, 2k-2, t) + t a(n-1, 2k, t). \tag{35}$$

This can be interpreted in the following way: Define the weight of an up-down path $P$ by $w(P) = t^{v(P)}$ where $v(P)$ is the number of down-steps $D_j : (i, j+1) \to (i, j)$ with odd $j$. The weight of a set of up-down-paths is the sum of the weights of its members.

Then $a(n,k,t)$ is the weight of all non-negative up-down-paths from $(0,0)$ to $(n,k)$.

This gives

**Proposition 1.**

$C_{k,n}(t)$ is the weight of all non-negative up-down-paths from $(0,0)$ to $(2n+k-1, k-1)$.

As an example, for $k=1$ the Narayana polynomials $C_n(t)$ can be interpreted as the weight of all up-down paths from $(0,0)$ to $(2n, 0)$.

If we denote an up-down path by the sequence of the $y$-coordinates of its steps we get

$w((1,0)) = 1 = C_1(t),$

$w((1,2,1,0)) + w((1,0,1,0)) = t + 1 = C_2(t),$

$w((1,2,3,2,1,0)) + w((1,2,1,2,1,0)) + w((1,0,1,0,1,0)) + w((1,0,1,2,1,0)) + w((1,2,1,0,1,0))$
$= t + t^2 + 1 + t + t = 1 + 3t + t^2 = C_3(t).$

As analog of (8) we get

**Proposition 2.**

$$\frac{1}{c(x,t)^{(k)}} + t^{\left\lfloor \frac{k+1}{2} \right\rfloor} x^k c(x,t)^{(k)} = h_k(x,t), \tag{36}$$

with

$$\begin{aligned} h_{2k}(x,t) &= L_k\left(1-(1+t)x, -tx^2\right), \\ h_{2k+1}(x,t) &= xtL_k\left(1-(1+t)x, -tx^2\right) + L_{k+1}\left(1-(1+t)x, -tx^2\right). \end{aligned} \tag{37}$$

**Remark**

For $t=1$ the first identity (37) reduces to $h_{2k}(x) = L_k\left(1-2x, -x^2\right) = L_{2k}(1,-x)$ by (13) and the second one to $L_{2k+2}(1,-x) + xL_{2k}(1,-x) = L_{2k+1}(1,-x) - xL_{2k}(1,-x) + xL_{2k}(1,-x) = L_{2k+1}(1,-x)$.

**Proof.**

For $k=1$ (36) reduces to .

$$\frac{1}{c_0(x,t)} + txc_0(x,t) = h_1(x,t) = 1 + (t-1)x. \tag{38}$$

This follows from (25) and (26) which give
$c_0(x,t) = 1 + xc_0(x,t)c_1(x,t) = 1 + xc_0(x,t)\left(1 - t + tc_0(x,t)\right)$ and therefore (38).

For $k=2$ we get

$$h_2(x,t) = 1 - (1+t)x. \tag{39}$$

The simplest proof uses a trick due to Max Alekseyev [1].

Since by (26) $\frac{1}{c_1(x,t)} + xtc_0(x,t) = 1$ and $\frac{1}{c_0(x,t)} + xc_1(x,t) = 1$ (38) implies

$$\left(\frac{1}{c_1(x,t)} + xtc_0(x,t)\right)\frac{1}{c_0(x,t)} + \left(\frac{1}{c_0(x,t)} + xc_1(x,t)\right)txc_0(x,t) = h_1(x,t) = 1 + (t-1)x,$$

which gives

$$h_2(x,t) = \frac{1}{c_0(x,t)c_1(x,t)} + tx^2 c_0(x,t)c_1(x,t) = 1 + (t-1)x - 2xt = 1 - (1+t)x.$$

Since $c(x,t)^{(2k)} = \left(c(x,t)^{(2)}\right)^k$ we get

$$h_{2k}(x,t) = \frac{1}{c(x,t)^{(2k)}} + t^k x^{2k} c(x,t)^{(2k)} = \frac{1}{\left(c(x,t)^{(2)}\right)^k} + t^k x^{2k} \left(c(x,t)^{(2)}\right)^k.$$

This implies as above that

$$h_{2k}(x,t) = L_k\left(1-(1+t)x, -tx^2\right). \tag{40}$$

To compute $h_{2k+1}(x,t)$ I use again Max Alekseyev's trick:

Let $a = \dfrac{1}{c_0(x,t)^{k+1} c_1(x,t)^k}$ and $b = t^{k+1} x^{2k+1} c_0(x,t)^{k+1} c_1(x,t)^k$.

Then

$$h_{2k+1}(x,t) = a+b = a\left(\frac{1}{c_1(x,t)} + xtc_0(x,t)\right) + b\left(\frac{1}{c_0(x,t)} + xc_1(x,t)\right)$$

$$= \left(c_0(x,t)tx + \frac{1}{c_1(x,t)}\right) \frac{1}{c_0(x,t)^{k+1} c_1(x,t)^k} + \left(c_1(x,t)x + \frac{1}{c_0(x,t)}\right) t^{k+1} x^{2k+1} c_0(x,t)^{k+1} c_1(x,t)^k$$

$$= \frac{tx}{c_0(x,t)^k c_1(x,t)^k} + \frac{1}{c_0(x,t)^{k+1} c_1(x,t)^{k+1}} + t^{k+1} x^{2k+2} c_0(x,t)^{k+1} c_1(x,t)^{k+1} + t^{k+1} x^{2k+1} c_0(x,t)^k c_1(x,t)^k$$

$$= txh_{2k}(x,t) + h_{2k+2}(x,t)$$

By (40) this gives

$$h_{2k+1}(x,t) = xtL_k\left(1-(1+t)x, -tx^2\right) + L_{k+1}\left(1-(1+t)x, -tx^2\right). \tag{41}$$

The polynomials $h_k(x,t)$ have degree $\left\lfloor \dfrac{k+1}{2} \right\rfloor$. The first polynomials $h_k(x,t)$ are

$h_1(x,t) = 1 - (1-t)x,$

$h_2(x,t) = 1 - (1+t)x,$

$h_3(x,t) = 1 - (2+t)x + (1-t)x^2,$

$h_4(x,t) = 1 - 2(1+t)x + (1+t^2)x^2,$

$h_5(x,t) = 1 - (3+2t)x + (3+t+t^2)x^2 - (1-t)x^3,$

$h_6(x,t) = 1 - 3(1+t)x + 3(1+t+t^2)x^2 - (1+t^3)x^3.$

Let

$$D_{K,M,t}(N) := \det\left(C_{K,i+j+M}(t)\right)_{i,j=0}^{N-1} \tag{42}$$

for $M \in \mathbb{Z}$ and integers $K, N \geq 1$. We set $D_{K,M,t}(0) = 1$ by definition.

**Theorem 3.**

*For $m \in \mathbb{N}$ and a positive integer $k$ we have*

$$D_{2k,1-k-m,t}(N) = 0 \quad \text{for } N = 1, 2, \cdots, m+k-1, \tag{43}$$

*and for all $n \in \mathbb{N}$*

$$D_{2k,1-k-m,t}(n+m+k) = (-1)^{\binom{m+k}{2}} t^{kn} D_{2k,1-k+m,t}(n). \tag{44}$$

For example, for $k=2$ and $m=1$ we get

$$\left(D_{4,-2,t}(n)\right)_{n\geq 0} = \left(1,0,0,-1,-t^2,t^4(1+t^2),t^8(1+t^2),-t^{12}(1+t^2+t^4),-t^{18}(1+t^2+t^4),t^{24}(1+t^2+t^4+t^6),\cdots\right)$$

$$\left(D_{4,0,t}(n)\right)_{n\geq 0} = \left(1,1,-(1+t^2),-t^2(1+t^2),t^4(1+t^2+t^4),t^8(1+t^2+t^4),-t^{12}(1+t^2+t^4+t^6),\cdots\right).$$

Note that by [4], Theorem 7.4. $D_{4,0,t}(2n) = (-1)^n t^{2(n^2-n)}\left(1+t^2+\cdots+t^{2n}\right)$ and
$D_{4,0,t}(2n+1) = (-1)^n t^{2n^2}\left(1+t^2+\cdots+t^{2n}\right).$

**Proof**

To prove (44) choose $s(x) = c(x,t)^{(2k)}$ and $M = m+k-1$, $N = n$. Then we get
$$\det\left(s_{i+j-M}\right)_{i,j=0}^{N+M} = \det\left(s_{i+j+1-m-k}\right)_{i,j=0}^{n+m+k-1} = D_{2k,1-k-m,t}(n+m+k).$$

By (36) we get

$t(x) = h_{2k}(x,t) - t^k x^{2k} c(x,t)^{(2k)}$ and therefore $t_{2k+r} = -t^k C_{2k,r}(t)$ and $t_r = 0$ for $k < r < 2k$.
This gives
$$\det\left(t_{i+j+M+2}\right)_{i,j=0}^{N-1} = \det\left(t_{i+j+m+k+1}\right)_{i,j=0}^{n-1} = \det\left(-t^k C_{2k,m-k+1}(t)\right)_{i,j=0}^{n-1} = (-1)^n t^{kn} D_{2k,m-k+1,t}(n).$$

For example, for $m=0$ we get

$$D_{2k,1-k,t}(kn) = (-1)^{n\binom{k}{2}} t^{k^2\binom{n}{2}}, \tag{45}$$
$$D_{2k,1-k,t}(n) = 0 \text{ else.}$$

**Theorem 4.**

*For positive integers $k$ and $m$ we have*
$$D_{2k-1,2-k-m,t}(N) = 0 \text{ for } N = 1,2,\cdots,m+k-2, \tag{46}$$

*and for all $n \in \mathbb{N}$*

$$D_{2k-1,2-k-m,t}(n+m+k-1) = (-1)^{\binom{m+k-1}{2}} t^{kn} D_{2k-1,1-k+m,t}(n). \tag{47}$$

For example, for $k=2$, $m=1$ we get

$$\left(D_{3,-1,t}(n)\right)_{n\geq 0} = \left(1,0,-1,-t^2,-t^4(-1+t),-t^8(-2+t)c-t^{12}(1-3t+t^2),-t^{18}(3-4t+t^2),\cdots\right),$$

$$\left(D_{3,0,t}(n)\right)_{n\geq 0} = \left(1,1,-1+t,t^2(-2+t),t^4(1-3t+t^2),t^8(3-4t+t^2),t^{12}(-1+6t-5t^2+t^3),\cdots\right).$$

Note that by [4], Theorem 5.2. $D_{3,0,t}(n) = \sum_{j=0}^{\lfloor \frac{n}{2} \rfloor} (-1)^j \binom{n-j}{j} t^{\binom{n}{2}-j}$.

**Proof**

Here we choose $s(x) = c(x,t)^{(2k-1)}$ and $M = m+k-2$. Then we get

$$\det(s_{i+j-M})_{i,j=0}^{N+M} = \det(s_{i+j+2-m-k})_{i,j=0}^{n+m+k-2} = D_{2k-1,2-k-m,t}(n+m+k-1).$$

On the other hand, we get $t(x) = h_{2k-1}(x,t) - t^k x^{2k-1} c(x,t)^{(2k-1)}$ and therefore

$t_n = -t^k s_{n+1-2k}$ for $n \geq k+1$.

This gives for $m \geq 1$

$$\det(t_{i+j+M+2})_{i,j=0}^{N-1} = \det(t_{i+j+m+k})_{i,j=0}^{n-1} = \det(-t^k C_{2k-1,m-k+1})_{i,j=0}^{n-1} = (-1)^n t^{kn} D_{2k-1,m-k+1,t}(n).$$

**Remark**

(33) and (45) for $k=1$ and $m=0$ and (47) for $m=1$ give

$D_{1,0,t}(n+1) = t^n D_{1,1,t}(n) = t^n D_{2,0,t}(n) = t^n t^{\binom{n}{2}} = t^{\binom{n+1}{2}}$. This gives the well-known analog of (19)

$$\det(C_{i+j}(t))_{i,j=0}^{n-1} = \det(C_{i+j+1}(t))_{i,j=0}^{n-1} = t^{\binom{n}{2}}. \tag{48}$$